\documentclass[12]{article}
\usepackage{fleqn}
\usepackage{graphicx}
\begin{document}
\Large
\begin{center}
{\bf A Classification of the Projective Lines over Small Rings\\
II. Non-Commutative Case}
\end{center}
\vspace*{.5cm}
\begin{center}
Metod Saniga,$^{1}$ Michel Planat$^{2}$ 
and Petr Pracna$^{3}$
\end{center}
\vspace*{.2cm} \normalsize
\begin{center}
$^{1}$Astronomical Institute, Slovak Academy of Sciences\\
SK-05960 Tatransk\' a Lomnica, Slovak Republic\\
(msaniga@astro.sk)

\vspace*{.4cm} $^{2}$Institut FEMTO-ST, CNRS, D\' epartement LPMO,
32 Avenue de
l'Observatoire\\ F-25044 Besan\c con Cedex, France\\
(planat@lpmo.edu)

\vspace*{.4cm} $^{3}$J. Heyrovsk\' y Institute of Physical
Chemistry, Academy of Sciences of the Czech Republic\\ Dolej\v
skova 3, CZ-182 23 Prague 8, Czech
Republic\\
(pracna@jh-inst.cas.cz)

\end{center}

\vspace*{.3cm} \noindent \hrulefill

\vspace*{.1cm} \noindent {\bf Abstract}

\noindent A list of different types of a projective line over non-commutative rings with unity of order up to thirty-one inclusive is given. Eight different types of such a line are found. With a single exception, the basic characteristics of the lines are identical to those of their commutative counterparts. The exceptional projective line is that defined over the non-commutative ring of order sixteen that features ten zero-divisors and it most pronouncedly differs from its commutative sibling in the number of shared points by the neighbourhoods of three pairwise distant points (three versus zero), that of ``Jacobson" points (zero versus five) and in the maximum number of mutually distant points (five versus three). 
\\



\noindent {\bf Keywords:} Projective Ring Lines -- Non-Commutative Rings of Small Orders

\vspace*{-.1cm} \noindent \hrulefill

\vspace*{.5cm}  \noindent
This is a short organic supplement to our recent paper [1], where the reader is referred to for the necessary background information on the concept of a projective ring line. Employing the definitions, symbols/notation and strategy of the above-mentioned paper, we have examined in detail the structure of different types of a projective line over small non-commutative rings with unity. As such rings are rather scarce for orders below thirty-two [2,3], we have found only eight line types, whose basic properties and all representative rings are listed in Table 1.  Here, the term projective line over $R$ means the {\it left}-line, i.e. the line whose points are given by the {\it left} equivalence classes ($\varrho a, \varrho b$),  where ($a, b$) is admissible over $R$ and $\varrho$ is a unit of $R$ [1]. 

\begin{table}[ht]
\begin{center}
\caption{The basic types of a projective line over small non-commutative rings with unity. The representative rings
are given in the notation of [2] (orders $\leq$16) and [3] (orders $>$16). For the line of 16/10 type, the numbers
given in brackets correspond to its commutative counterpart.} 
\vspace*{0.3cm}
{\begin{tabular}{|||l|c|c|c|c|c|c|c|l|||} \hline \hline \hline
\multicolumn{1}{|||c|}{} & \multicolumn{7}{|c|}{} &  \multicolumn{1}{|c|||}{}\\
\multicolumn{1}{|||c|}{Line} & \multicolumn{7}{|c|}{Cardinalities
of Points} &
\multicolumn{1}{|c|||}{Representative} \\
\multicolumn{1}{|||c|}{Type} & \multicolumn{7}{|c|}{} &  \multicolumn{1}{|c|||}{Rings} \\
\cline{2-8}
 &Tot& TpI & 1N & $\cap$2N & $\cap$3N & Jcb & MD &   \\
 \hline \hline
27/15 & 48 & 42 & 20 & 6 & 0 & 2 & 4 & 3.20  \\
\hline
\hline 
24/20 & 72 & 44 & 47 & 28 & 12 & 3 & 3 & 3.22  \\
\hline
\hline 
16/4 & 20 & 20 & 3 & 0 & 0 & 3 & 5 & 5.105 \\
\hline
16/8 & 24 & 24 & 7 & 0 & 0 & 7 & 3 & 4.83,\,4.117,\,5.98,\,5.101   \\
16/10 & 35(30) & 26 & 18(13) & 9(4) & 3(0) & 0(5) & 5(3) & 5.96 \\
16/12 & 36 & 28 & 19 & 8 & 0 & 3 & 3 & 4.68,\,4.73,\,5.97,\,5.100,\,5.104,\,5.106  \\
16/14 & 54 & 30 & 37 & 24 & 12 & 1 & 3 & 5.113  \\
\hline
\hline 
\hspace*{0.05cm} 8/6 & 18 & 14 & 9 & 4 & 0 & 1 & 3 & 3.11 \\
\hline \hline \hline
\end{tabular}}
\end{center}
\end{table}

Comparing Table 1 with Table 3 of [1], one readily notices the exceptional character of the line of 16/10 type, which differs from its commutative counterpart in a number of aspects. This difference stems from the fact that its base ring [3; pp.\,433, 531], although having no two-sided ideals, has {\it three} proper maximal right- (and also left-) ideals to be compared with only {\it two} proper (and, of course, two-sided) maximal ideals of the corresponding commutative rings ($\cong GF(4)\otimes Z_{4}$ or $GF(4)\otimes GF(2)[x]/\langle x^{2}\rangle$) [1]. A deeper insight into this intriguing difference is acquired when we try to pass to the {\it right}-line over $R$, i.e., the line whose points are regarded as right equivalence classes ($a \varrho, b \varrho$). For all rings under consideration except the 16/10 one, this right-line was found to exist and possess the same properties as its left-companion. In the 16/10 case, however, this concept breaks down due to the fact that equivalence classes are not of the same cardinality, which renders it impossible to define consistently the notions of neighbour/distant [4]. 

Attacking the next order, thirty-two, in the hierarchy seems to represent a truly big computational challenge; this not only because of a large number (ninety-nine in total) of distinct non-commutative rings with unity there [5], but also, and mainly, due to the fact that their addition and multiplication tables have not been published/available yet. 

\vspace*{.7cm} \noindent \Large {\bf Acknowledgements}
\normalsize

\vspace*{.3cm} \noindent This work was partially supported by the
Science and Technology Assistance Agency under the contract $\#$
APVT--51--012704, the VEGA project $\#$ 2/6070/26 (both from
Slovak Republic), the trans-national ECO-NET project $\#$
12651NJ ``Geometries Over Finite Rings and the Properties of
Mutually Unbiased Bases" (France) and by the project 1ET400400410 of 
the Academy of Sciences of the Czech Republic.

\end{document}